# Multiplicative orders of elements in Conway's towers of finite fields

Roman Popovych

Department of Computer Technologies, Lviv Polytechnic National University, Bandery Str. 12, Lviv 79013, Ukraine
e-mail: rombp07@gmail.com

**Abstract.** We give a lower bound on multiplicative orders of some elements in defined by Conway towers of finite fields of characteristic two and also formulate a condition under that these elements are primitive.

**Classification code:** AMS 11T30

**Keywords:** finite field; multiplicative order, Conway's tower

Elements with high multiplicative order are often needed in several applications that use finite fields [11]. Ideally we want to have a possibility to obtain a primitive element for any finite field. However, if we have no the factorization of the order of finite field multiplicative group, it is not known how to rich the goal. That is why one considers less ambitious question: to find an element with provable high order. It is sufficient in this case to obtain a lower bound on the order. The problem is considered both for general and special finite fields [1,3,7, 12,13].

Another less ambitious, but supposedly more important question is to find primitive elements for a class of special finite fields. A polynomial algorithm that find a primitive element in finite field of small characteristic is described in [8]. However, the algorithm relies on two unproved assumptions and is not supported by any computational example. Our paper can be considered as a step towards this direction. We give a lower bound on multiplicative order of some elements in binary recursive extensions of finite fields defined by Conway [4,5,14] and also formulate a condition under that these elements are primitive. $F_q$ denotes finite field with $q$ elements. More precisely, we consider the following finite fields:

$c_{-1} = 1, L_{-1} = F_2(c_{-1}) = F_2$;

for $i \geq -1$, $L_{i+1} = L_i(c_{i+1})$, where $c_{i+1}$ satisfies the equation $c_{i+1}^2 + c_{i+1} + \prod_{j=-1}^{i} c_j = 0$.

So, the following tower of finite fields of characteristic two arises:

$$L_{-1} = F_2(c_{-1}) = F_2 \subset L_0 = F_2(c_0) \subset L_1 = L_0(c_1) \subset L_2 = L_1(c_2) \subset ...$$

Such construction is very attractive from point of view of applications, since one can perform operations with finite field elements recursively, and therefore effectively [9].

It is easy to verify directly the following facts: element $c_0$ is primitive in $L_0$, and element $c_1$ is primitive in $L_1$. On the other hand, H. Lenstra [10, exercise 2] showed: if $i \geq 2$, then element $c_i$ is not primitive in $L_i$. Some primitive elements for the fields $L_2$, $L_3$, $L_4$ are found in [2]. Therefore, for $i \geq 2$, the following questions arise: 1) give a lower bound on the multiplicative order $O(c_i)$ of element $c_i$; 2) what elements are primitive in the fields $K_i$. We answer partially the questions in theorems 3, 4 and corollaries 2, 3, 4, 5.

---

Date: 07 September 2015

**1. Preliminaries**

Point out that, for $i \geq 0$, the number of elements of the multiplicative group $L_i^* = L_i \setminus \{0\}$ equals to $2^{2^{i+1}} - 1$. If to denote the Fermat numbers $N_j = 2^{2^j} + 1$ ($j \geq 0$), then the cardinality of $L_i^*$ equals to $\prod_{j=0}^{i} N_j$. We use for $k \geq 0$ the denotation $a_k = \prod_{j=0}^{k} c_j$. $(a,b)$ denotes the greatest common divisor of the numbers $a$ and $b$. $a \mid b$ means that $a$ divides $b$.

**Lemma 1** [6, section 1.3]. *For $i \geq 1$, $N_i = \prod_{j=0}^{i-1} N_j + 2$.*

**Lemma 2** [6, section 1.3]. *Any two Fermat numbers are coprime.*

**Lemma 3.** *For $j \geq 2$, a divisor $\alpha > 1$ of the number $N_j$ is of the form $\alpha = l \cdot 2^{j+2} + 1$, where $l$ is a positive integer.*

**Proof.** The result obtained by Euler and Lucas (see [6, theorem 1.3.5]) states: for $j \geq 2$, a prime divisor of the number $N_j$ is of the form $l \cdot 2^{j+2} + 1$, where $l$ is a positive integer. Clearly a product of two numbers of the specified form is a number of the same form. Hence, the result follows.

**Lemma 4.** *For $i \geq 2$, $(N_i + 1, N_j) = 1$ for $1 \leq j \leq i-1$.*

**Proof.** According to lemma 1, $N_i + 1 = \prod_{j=0}^{i-1} N_j + 3$. Common divisor of numbers $N_i + 1$ and $\prod_{j=0}^{i-1} N_j$ divides their difference that equals to 3. As $N_0 = 3$, we have $(N_i + 1, \prod_{j=0}^{i-1} N_j) = 3$. Since, according to lemma 2, numbers $N_j$ are coprime, $(N_i + 1, N_j) = 1$ for $i \geq 2$ and $1 \leq j \leq i-1$.

**Lemma 5.** *For $i \geq 1$, the equalities*
$$(c_i)^{N_i} = a_{i-1} \quad (1)$$
*and*
$$(a_i)^{N_i} = (a_{i-1})^{N_i+1}. \quad (2)$$
*are true.*

**Proof.** Show first that (1) holds. Indeed, note that $c_i$ is a root of the equation $x^2 + x + a_{i-1} = 0$ over the field $L_{i-1}$. One can verify directly, that $c_i + 1$ is also a root of this equation. Then $c_i$ and $c_i + 1$ are conjugates [11] over $L_{i-1} = F_{2^{2^i}}$, that is $(c_i)^{2^{2^i}} = c_i + 1$. Therefore, $(c_i)^{2^{2^i}+1} = (c_i+1)c_i = a_{i-1}$, and the equality (1) is true. Applying (1), we obtain $(a_i)^{N_i} = (c_i a_{i-1})^{N_i} = (a_{i-1})^{N_i+1}$. Hence, (2) is true as well.

If $u_j$ is a sequence of integers and $s > t$, then we consider below the empty product $\prod_{j=s}^{t} u_j = 1$.

**Lemma 6.** *For $k \geq 0$, the following equalities are true*

$$(c_i)^{\prod_{j=0}^{k} N_{i-j}} = (a_{i-k-1})^{\prod_{j=1}^{k}(N_{i-j}+1)} \tag{3}$$

*and*

$$(a_i)^{\prod_{j=0}^{k} N_{i-j}} = (a_{i-k-1})^{\prod_{j=0}^{k}(N_{i-j}+1)} \tag{4}$$

*for $i > k$.*

**Proof.** We will prove this by induction on $k$. For $k = 0$ (and for $i \geq 1$) we have the equalities (1) and (2).

Assume the equalities (3) and (4) hold for $k-1$, that is

$$(c_i)^{\prod_{j=0}^{k-1} N_{i-j}} = (a_{i-(k-1)-1})^{\prod_{j=1}^{k-1}(N_{i-j}+1)} \tag{5}$$

and

$$(a_i)^{\prod_{j=0}^{k-1} N_{i-j}} = (a_{i-(k-1)-1})^{\prod_{j=0}^{k-1}(N_{i-j}+1)}. \tag{6}$$

Show first that (3) is also true for $k$. Indeed, applying (5) and (2), we obtain

$$(c_i)^{\prod_{j=0}^{k} N_{i-j}} = \left((c_i)^{\prod_{j=0}^{k-1} N_{i-j}}\right)^{N_{i-k}} = \left((a_{i-k})^{N_{i-k}}\right)^{\prod_{j=1}^{k-1}(N_{i-j}+1)} = (a_{i-k-1})^{\prod_{j=1}^{k}(N_{i-j}+1)}.$$

Show now that (4) is also true for $k$. Really, applying (6) and (2), we obtain

$$(a_i)^{\prod_{j=0}^{k} N_{i-k}} = \left((a_i)^{\prod_{j=0}^{k-1} N_{i-k}}\right)^{N_{i-k}} = \left((a_{i-k})^{N_{i-k}}\right)^{\prod_{j=0}^{k-1}(N_{i-j}+1)} = (a_{i-k-1})^{\prod_{j=0}^{k}(N_{i-j}+1)}.$$

**Lemma 7.** *Let $K \subset L$ be a tower of fields. Let $x \in L \setminus K$ and $m$ be the smallest positive integer, satisfying the condition $x^m \in K$. If $x^n \in K$ for a positive integer $n$, then $m \mid n$.*

**Proof.** One can write $n = um + v$, where $0 \leq v < m$. Then $x^n = (x^m)^u \cdot x^v$, and, therefore, $x^v \in K$. As $m$ is the smallest positive integer with the condition $x^m \in K$ and $v < m$, we have $v = 0$, and the result follows.

**Lemma 8.** *For $u \geq 1$, if $(c_u)^l \in L_{u-1}$, where $l$ is a positive integer, then $(l, N_u) > 1$.*

**Proof.** The equality (1) gives $(c_u)^{N_u} = a_{u-1} \in L_{u-1}$. According to lemma 7, if $d$ is the smallest positive integer with $(c_u)^d \in L_{u-1}$, then $d \mid N_u$ and $d \mid l$. Clearly $d > 1$. Therefore $(l, N_u) \geq d > 1$.

**Lemma 9.** *Let $L_1 \subset L_2$ be a tower of fields and $b \in L_2^*$. Let $b^r = a \in L_1^*$ and $r$ be the smallest positive integer with $b^r \in L_1^*$. Then $O(b) = r \cdot O(a)$.*

**Proof.** As $b^{O(b)} = 1 \in L_1^*$, one has $O(b) \geq r$. Write $O(b) = sr + t$, where $s \in N$ and $0 \leq t < r$. Then
$$1 = b^{O(b)} = b^{sr+t} = a^s b^t.$$
Hence, $b^t = a^{-s} \in L_1^*$. By definition of $r$, it is possible only for $t = 0$. We have $a^s = 1$, $s \geq O(a)$ and $O(b) = sr \geq r \cdot O(a)$. From the other side, $b^{r \cdot O(a)} = a^{O(a)} = 1$, and therefore $O(b) = r \cdot O(a)$.

**Theorem 1.** *For $i \geq 2$, $(c_i)^{\prod_{j=0}^{k} N_{i-j}} \in L_{i-k-1} \setminus L_{i-k-2}$ for $0 \leq k \leq i-1$.*

**Proof.** Applying equality (3), one obtains
$$(c_i)^{\prod_{j=0}^{k} N_{i-j}} = (c_{i-k-1})^{\prod_{j=1}^{k} (N_{i-j}+1)} (a_{i-k-2})^{\prod_{j=1}^{k} (N_{i-j}+1)}. \tag{7}$$

Obviously $(c_{i-k-1})^{\prod_{j=1}^{k}(N_{i-j}+1)} \in L_{i-k-1}$ and $(a_{i-k-2})^{\prod_{j=1}^{k}(N_{i-j}+1)} \in L_{i-k-2}$. Hence, the product on the right hand of equality (7) belongs to $L_{i-k-1}$. For $1 \leq j \leq k$, according to lemma 4, $(N_{i-j}+1, N_{i-k-1}) = 1$ and $(\prod_{j=1}^{k}(N_{i-j}+1), N_{i-k-1}) = 1$. Then, by lemma 8, $(c_{i-k-1})^{\prod_{j=1}^{k}(N_{i-j}+1)} \notin L_{i-k-2}$. Therefore, the element $(c_{i-k-1})^{\prod_{j=1}^{k}(N_{i-j}+1)} (a_{i-k-2})^{\prod_{j=1}^{k}(N_{i-j}+1)}$ does not belong to $L_{i-k-2}$.

**Theorem 2.** *For $i \geq 2$, $(a_i)^{\prod_{j=0}^{k} N_{i-j}} \in L_{i-k-1} \setminus L_{i-k-2}$ for $0 \leq k \leq i-1$.*

**Proof.** Applying equality (4), we have
$$(a_i)^{\prod_{j=0}^{k} N_{i-j}} = (c_{i-k-1})^{\prod_{j=0}^{k} (N_{i-j}+1)} (a_{i-k-2})^{\prod_{j=0}^{k} (N_{i-j}+1)}. \tag{8}$$

Obviously $(c_{i-k-1})^{\prod_{j=0}^{k}(N_{i-j}+1)} \in L_{i-k-1}$ and $(a_{i-k-2})^{\prod_{j=0}^{k}(N_{i-j}+1)} \in L_{i-k-2}$. Hence, the product on the right hand of the equality (8) belongs to $L_{i-k-1}$. For $0 \leq j \leq k$, according to lemma 4, $(N_{i-j}+1, N_{i-k-1}) = 1$ and $(\prod_{j=0}^{k}(N_{i-j}+1), N_{i-k-1}) = 1$. So, by lemma 8, $(c_{i-k-1})^{\prod_{j=0}^{k}(N_{i-j}+1)} \notin L_{i-k-2}$. Hence, the element $(c_{i-k-1})^{\prod_{j=0}^{k}(N_{i-j}+1)} (a_{i-k-2})^{\prod_{j=0}^{k}(N_{i-j}+1)}$ does not belong to $L_{i-k-2}$.

## 2. Lower bound on multiplicative orders of elements

We give in this section in corollary 2 a lower bound on multiplicative order of elements $c_i$, $a_i$ and also formulate in corollary 3 a condition under that these elements are primitive.

**Theorem 3.** *For $i \geq 2$, the following statements hold:*

$$(a)\ O(c_i) = \prod_{j=1}^{i} \alpha_j,\ \text{where}\ \alpha_j \mid N_j,\ \alpha_j > 1;$$

$$(b)\ O(a_i) = \prod_{j=1}^{i} \beta_j,\ \text{where}\ \beta_j \mid N_j,\ \beta_j > 1.$$

**Proof.** (a) Define recursively the sequence $\alpha_i, \ldots, \alpha_{1i}$ of positive integers as follows. $\alpha_i$ is the smallest integer satisfying the condition $(c_i)^{\alpha_i} \in L_{i-1}$. If $\alpha_i, \ldots, \alpha_{i-j}$, where $0 \leq j \leq i-2$, are already known, then $\alpha_{i-j-1}$ is the smallest positive integer such that the condition

$$\{(c_i)^{\prod_{k=i-j}^{i} \alpha_k}\}^{\alpha_{i-j-1}} \in L_{i-j-2}$$

holds.

Since the cardinality of the group $L_i^*$ equals to $\prod_{j=0}^{i} N_j$ and the cardinality of the group $L_{i-1}^*$ equals to $\prod_{j=0}^{i-1} N_j$, we have that the number of elements of the factor-group $L_i^* / L_{i-1}^*$ equals to $N_i$. If $d$ is the coset of $c_i$ in the factor-group, then $\alpha_i = O(d)$ and as a consequence of Lagrange's theorem for finite groups $\alpha_i \mid N_i$. Clearly, $\alpha_i > 1$. $(c_i)^{\alpha_i} \in L_{i-1} \setminus L_{i-2}$, because according to theorem 2 $(c_i)^{N_i} \in L_{i-1} \setminus L_{i-2}$. Indeed, if to suppose $(c_i)^{\alpha_i} \in L_{i-2}$, then $[(c_i)^{\alpha_i}]^{N_i/\alpha_i} = (c_i)^{N_i} \in L_{i-2}$ - a contradiction. According to lemma 9, $O(c_i) = \alpha_i O((c_i)^{\alpha_i})$.

Analogously to the previous, one can show in the case of finding $\alpha_{i-j-1}$ that $\alpha_{i-j-1} \mid N_{i-j-1}$ ($\alpha_{i-j-1} > 1$) and $\{(c_i)^{\prod_{k=i-j}^{i} \alpha_{i-k}}\}^{\alpha_{i-j-1}} \in L_{i-j-2} \setminus L_{i-j-3}$. According to lemma 9, $O((c_i)^{\alpha_i \ldots \alpha_{i-j}}) = \alpha_{i-j-1} O((c_i)^{\alpha_i \ldots \alpha_{i-j} \alpha_{i-j-1}})$.

Since, according to (3),

$$(c_i)^{\prod_{j=0}^{i-1} N_{i-j}} = ((a_0)^{N_1+1})^{\prod_{j=1}^{i-2}(N_{i-j}+1)} = 1,$$

one obtains $O(c_i) \mid \prod_{j=0}^{i-1} N_{i-j}$, and $O(c_i) = \alpha_i \ldots \alpha_1$.

(b) The proof is analogues to the previous one, using theorem 2 instead of theorem 1.

**Corollary 1.** *For* $i \geq 2$, $O(c_i c_0) = N_0 O(c_i)$ *and* $O(a_i a_0) = N_0 O(a_i)$.

**Proof.** Note that $O(c_0) = N_0$. Since, according to theorem 3, $O(c_i)$ divides $\prod_{j=1}^{i} N_j$, and lemma 2 gives $(\prod_{j=1}^{i} N_j, N_0) = 1$, we have $(O(c_i), O(c_0)) = 1$. Then $O(c_i c_0) = O(c_i) O(c_0)$, and the result for $c_i c_0$ follows. The proof for $a_i a_0 = a_i c_0$ is analogous.

**Corollary 2.** *The multiplicative order of the element* $c_i$ *and the element* $a_i$ *equals to* $\prod_{j=1}^{i} N_j$ *for* $2 \leq i \leq 4$ *and is at least* $\prod_{j=1}^{4} N_j \cdot \prod_{j=5}^{i} (2^{j+2} + 1)$ *for* $i \geq 5$.

**Proof.** Consider the formulas for the multiplicative order of $c_i$ and $a_i$ given in theorem 3. It is known that for $1 \leq j \leq 4$ the Fermat numbers $N_1 = 5$, $N_2 = 17$, $N_3 = 257$, $N_4 = 65537$ are prime [6, table 1.3]. Therefore, $\alpha_j = \beta_j = N_j$ for $1 \leq j \leq 4$. According to lemma 3, for $j \geq 5$, $\alpha_j, \beta_j \geq 2^{j+2} + 1$.

**Theorem 4.** *Let* $i \geq 5$. *For* $5 \leq j \leq i$, *if* $\alpha_j = N_j$ *is the smallest positive integer satisfying the condition* $(c_j)^{\alpha_j} \in L_{j-1}$, *then* $O(a_i) = O(c_i) = \prod_{j=1}^{i} N_j$.

**Proof.** For element $a_i$ the proof is by induction on $i \geq 5$. Note that $\alpha_j$ is the smallest positive integer with $(c_j)^{\alpha_j} \in L_{j-1}$ iff $\alpha_j$ is the smallest positive integer with $(a_j)^{\alpha_j} \in L_{j-1}$.

For $i = 5$ we have from equality (2) $(a_5)^{N_5} = (a_4)^{N_5+1}$. Then, according to lemma 9, $O(a_5) = N_5 O((a_4)^{N_5+1})$. We have $O(a_4) = \prod_{j=1}^{4} N_j$ by corollary 2 and $(N_5 + 1, \prod_{j=1}^{4} N_j) = 1$ by lemma 4. Now use the well known fact that raising an element of a group to a power relatively prime to its order does not change the order. One obtains $O((a_4)^{N_5+1}) = O(a_4)$ and $O(a_5) = \prod_{j=1}^{5} N_j$.

Assume the statement of the theorem is true for $i-1$. For $i$ we have from equality (2) $(a_i)^{N_i} = (a_{i-1})^{N_i+1}$. Then, according to lemma 9, $O(a_i) = N_i O((a_{i-1})^{N_i+1})$. As $O(a_{i-1}) = \prod_{j=1}^{i-1} N_j$ by induction assumption and $(N_i + 1, \prod_{j=1}^{i-1} N_j) = 1$ by lemma 4, one has, analogously to the previous, $O((a_{i-1})^{N_i+1}) = O(a_{i-1})$ and $O(a_i) = \prod_{j=1}^{i} N_j$.

To complete the proof note that, according to equality (1) and lemma 9, $O(c_i) = N_i O(a_{i-1}) = O(a_i)$.

Remark that, if the condition of theorem 4 is true, then one has the following chain of cyclic subgroups:
$$\langle c_i \rangle = \langle a_i \rangle \supset \langle c_{i-1} \rangle = \langle a_{i-1} \rangle \supset ... \supset \langle c_2 \rangle = \langle a_2 \rangle \supset \langle a_1 \rangle.$$
At the same time $\langle a_1 \rangle \neq \langle c_1 \rangle$, because $O(c_1) = 15$, $O(a_1) = O(c_1 c_0) = 5$.

As a consequence of theorem 4 and corollary 1, one has the following corollary.

**Corollary 3.** Let $i \geq 5$. For $5 \leq j \leq i$, if $\alpha_j = N_j$ is the smallest positive integer with $(c_j)^{\alpha_j} \in L_{j-1}$, then $c_i c_0$ and $a_i a_0$ is primitive.

**Proof.** Since $O(c_i c_0) = O(a_i a_0) = \prod_{j=0}^{i} N_j$, the result follows.

**Theorem 5.** *For $5 \leq j \leq 11$, the number $\alpha_j = N_j$ is the smallest positive integer with $(c_j)^{\alpha_j} \in L_{j-1}$.*

**Proof.** For $5 \leq j \leq 11$, the Fermat numbers are completely factored into primes [6]. Correspondent factorizations and connected with them results of verification are given below. Note that to prove the fact: $N_j$ is the smallest positive integer with $(c_j)^{\alpha_j} \in L_{j-1}$, it is enough to verify $c_j^{N_j/p} \notin L_{j-1}$ for any prime divisor $p$ of $N_j$. Really, if element $c_j$ in the power $N_j/p$ does not belong to $L_{j-1}$, then element $c_j$ in the power of any divisor $N_j/(pq)$ of $N_j/p$ does not belong to $L_{j-1}$ as well.

For $j = 5$, one has the factorization $N_5 = 641 \cdot 6700417$. According to equality (1), $(c_5)^{N_5} = c_4 c_3 c_2 c_1 c_0 \in L_4$. We have verified that $\alpha_5 = N_5$ is the smallest positive integer with $(c_5)^{\alpha_5} \in L_4$.

For $j = 6$, one has the factorization $N_6 = 274177 \cdot 67280421310721$. According to equality (1), $(c_6)^{N_6} = c_5 c_4 c_3 c_2 c_1 c_0 \in L_5$. We have verified that $\alpha_6 = N_6$ is the smallest positive integer with $(c_6)^{\alpha_6} \in L_5$.

For $j = 7$ one has the factorization $N_7 = 59649589127497217 \cdot 5704689200685129054721$. According to equality (1), $(c_7)^{N_7} = c_6 c_5 c_4 c_3 c_2 c_1 c_0 \in L_6$. We have verified that $\alpha_7 = N_7$ is the smallest positive integer with $(c_7)^{\alpha_7} \in L_6$.

For $j = 8$, one has the factorization $N_8 = 1238926361552897 \cdot P_{62}$, where $P_{62}$ is a prime with 62 decimal digits,
$P_{62}$=93461639715357977769163558199606896584051237541638188580280321.
According to equality (1), $(c_8)^{N_8} = c_7 c_6 c_5 c_4 c_3 c_2 c_1 c_0 \in L_7$. We have verified that $\alpha_8 = N_8$ is the smallest positive integer with $(c_8)^{\alpha_8} \in L_7$.

For $j = 9$, one has the factorization
$N_9 = 2424833 \cdot 7455602825647884208337395736200454918783366342657 \cdot P_{99}$, where $P_{99}$ is a prime number with 99 decimal digits,

P$_{99}$=741640062627530801524787141901937474059940781097519023905821316144415759504
70500809281871169394073 7.

According to equality (1), $(c_9)^{N_9} = c_8 c_7 c_6 c_5 c_4 c_3 c_2 c_1 c_0 \in L_8$. We have verified that $\alpha_9 = N_9$ is the smallest positive integer with $(c_9)^{\alpha_9} \in L_8$.

For $j = 10$, one has the factorization

$N_{10} = 45592577 \cdot 6487031809 \cdot 4659775785220018543264560743076778192897 \cdot P_{252}$, where P$_{252}$ is a prime with 252 decimal digits,

P$_{252}$=130439874405488189727484768796509903946608530841611892186895295776832416251471863574140227977573104895898783928842923844831149032913798729088601617946094119449010595906710130531906171018354491609619193912488538116080712299672322806217820753127014424577.

According to equality (1), $(c_{10})^{N_{10}} = c_9 c_8 c_7 c_6 c_5 c_4 c_3 c_2 c_1 c_0 \in L_9$. We have verified that $\alpha_{10} = N_{10}$ is the smallest positive integer with $(c_{10})^{\alpha_{10}} \in L_9$.

For $j = 11$, one has the factorization

$N_{11} = 319489 \cdot 974849 \cdot 167988556341760475137 \cdot 3560841906445833920513 \cdot P_{564}$, where P$_{564}$ is a prime with 564 decimal digits,

P$_{564}$=17346244717914755543025897086430977837742184472366408464934701906136357919287910885759103833040883717798381086845154642194071297830613418986428082601454275870858924387368556397311894886939915854550661147420216132557017260564139394366945793220968665108959685482705388072645828554151936401912464931182546092879815733057795573358504982279280090942872567591518912118622751714319229788100979251036035496917279912663527358783236647193154777091427745377038294584918917590325110939381322486044298573971650711059244462177542540706913047034664643603491382441723306598834177.

According to equality (1), $(c_{11})^{N_{11}} = c_{10} c_9 c_8 c_7 c_6 c_5 c_4 c_3 c_2 c_1 c_0 \in L_{10}$. We have verified that $\alpha_{11} = N_{11}$ is the smallest positive integer with $(c_{11})^{\alpha_{11}} \in L_{10}$.

**Corollary 4.** *For $2 \leq i \leq 11$, the multiplicative order of element $c_i$ and element $a_i$ equals to $\prod_{j=1}^{i} N_j$.*

**Proof.** The case $2 \leq i \leq 4$ follows from corollary 2. The case $5 \leq i \leq 11$ follows from theorem 4 and theorem 5.

As a consequence of corollary 1 and corollary 4, one has the following corollary.

**Corollary 5.** For $2 \leq i \leq 11$, the element $c_i c_0$ and the element $a_i a_0$ is primitive.

Consider at the end of the paper the following example concerning with the multiplicative group of the field $K_2$. We have $O(c_2) = 5 \cdot 17 = 85$. The element $c_2 c_0$ is primitive in the field, i.e. $O(c_2 c_0) = 3 \cdot 5 \cdot 17 = 255$. Take also the element $c_2 + c_1 + 1$. Since $c_2 + c_1 + 1 = (c_2)^5$, we obtain $O(c_2 + (c_1 + 1)) = 17$.